\documentclass[a4paper, 11pt]{amsart}
\usepackage[utf8]{inputenc}
\usepackage{amsfonts, amsthm, amsmath, amssymb,bm}

\usepackage[T1]{fontenc}
\usepackage{verbatim}
\usepackage{amssymb,amscd,amsmath,amsthm}
\usepackage{mathtools}
\usepackage{enumerate}
\usepackage[subnum]{cases}

\usepackage{fancyhdr}

\usepackage{bbm}

\theoremstyle{plain}
\newtheorem{thm}{Theorem}[section]
\newtheorem{lem}[thm]{Lemma}

\newtheorem{cor}[thm]{Corollary}

\theoremstyle{definition}

\usepackage{mdframed}

\newcommand{\Z}{\mathbb{Z}}

\newcommand{\rmod}{\!\!\!\!\pmod}

\usepackage[square,sort,comma,numbers]{natbib}
\usepackage{graphicx}

\title[On the least squarefree number in an arithmetic progression]{A note on the least squarefree number in an arithmetic progression}

\author{Ramon M. Nunes} \thanks{This work is supported by the DFG-SNF lead agency program grant 200021L-153647.}

\address{EPFL SB MATHGEOM TAN\\ Station 8\\ CH-1015 Lausanne\\ Switzerland}
\email{ramon.moreiranunes@epfl.ch}

\usepackage{hyperref}

\begin{document}

\maketitle

\begin{abstract}
We prove an asymptotic formula for squarefree in arithmetic progressions with squarefree moduli, improving previous results by Prachar. The main tool is  an estimate for counting solutions of a congruence inside a box that goes beyond what can be obtained by using the Weil bound.
\end{abstract}

\section{Introduction}

Let $\mu$ denote the M\"obius function, \textit{i.e.}~$\mu$ is the multiplicative function such that for every prime number $p$ and every positive integer $j$, one has,
$$
\mu(p^{j})=
\begin{cases}
-1,\,\text{if }j=1,\\
\;\;\;0,\,\text{otherwise}.
\end{cases}
$$
We remark that $\mu^2(n)=1$ if $n$ is squarefree and $\mu^2(n)=0$ otherwise. In this paper we are concerned with the distribution of squarefree numbers in arithmetic progressions. By the above discussion, this is equivalent to studying the distribution of the $\mu^2$ function in arithmetic progressions.

In this direction, a result of Prachar \citep{prachar1958kleinste}, subsequently improved by Hooley \citep{hooley1975note} says that
\begin{equation}\label{Hoo-sqf}
\sum_{\substack{n\leq x\\n \equiv a\rmod q}}\mu^2(n) = \frac{1}{\varphi(q)}\sum_{\substack{n\leq x\\(n,q)=1}}\mu^2(n) + O\left(\frac{X^{1/2}}{q^{1/2}} + q^{1/2+\epsilon}\right).
\end{equation}
Here and throughout the article, $\epsilon$ denotes a small constant that might vary from line to line and the implied constants in the symbols $O$ and $\ll$ are allowed to depend on $\epsilon$.

It follows from \eqref{Hoo-sqf} that the sequence of squarefree numbers $\leq X$ is well distributed in arithmetic progressions modulo $q$ whenever
\begin{equation}\label{q<X23}
q\leq X^{2/3-\epsilon}.
\end{equation}
Even though it is largely believed that one should be able to replace $2/3$ by $1$ in the above inequality, this constant has resisted any improvement until very recently.

The author \citep{nunes2016squarefree} proved, using more sophisticated techniques than those contained here, that if one restricts to prime values of $q$, the exponent in \eqref{q<X23} can be improved to $13/19$. In the present paper we show how to further improve this constant and at the same time relax the condition on $q$. Our main result is the following:
\begin{thm}\label{25/36}
Let $\epsilon>0$. Then there exists $\delta=\delta(\epsilon)>0$ such that, uniformly for $X\geq 2$, integers $a$ and squarefree numbers $q$ coprime with $a$ satisfying
$$
q\leq X^{\frac{25}{36}-\epsilon},
$$
we have
\begin{equation}\label{asform}
\sum_{\substack{n\leq X\\n \equiv a\rmod q}}\mu^2(n) = \frac{1}{\varphi(q)}\sum_{\substack{n\leq X\\(n,q)=1}}\mu^2(n) + O\left(\frac{X^{1-\delta}}{q}\right).
\end{equation}
In other terms, the value $\Theta=\frac{25}{36}$ is an exponent of distribution for the characteristic function of the sequence of squarefree numbers $\mu^2$ restricted to squarefree moduli.
\end{thm}
Alternatively, one can ask the simpler question of when is the left-hand side of \eqref{Hoo-sqf} nonzero. This is equivalent to study the least squarefree number in an arithmetic progression. We let $n(q,a)$ denote the least positive squarefree number which is congruent to $a$ modulo $q$. Prachar's result implies
\begin{equation}\label{3/2}
n(q,a)\ll q^{\frac32+\epsilon}.
\end{equation}
This result was improved by Erd\"os \citep{erdos1960kleinste}, who essentially proved that $n(q,a)=o(q^{3/2})$ and then later by Heath-Brown \citep{heath1982least}, who showed the upper bound
\begin{equation}\label{13/9}
n(q,a)\ll q^{\frac{13}{9}+\epsilon}.
\end{equation}
It is a direct consequence of Theorem \ref{25/36} that we can also improve this inequality for squarefree values of $q$. Indeed we have the following:


\begin{cor}\label{least}
For every $\epsilon>0$, we have the inequality
$$
n(q,a)\ll q^{\frac{36}{25} +\epsilon},
$$
uniformly for $q$ squarefree and $a$ coprime with $q$.
\end{cor}

The key input comes from estimates for the number of solutions to a congruence inside a dyadic box that follow from the work of Pierce \citep{pierce2005part}.

Let $q$ be a positive integer and let $a\in (\Z/q\Z)^{\ast}$. Let $u>0$ and $v$ be nonzero integers and let $M$ and $N$ be real numbers such that $M, N \geq 1$. We consider the counting function
$$
S_{u,v}(M,N,q,a) := \#\{m\leq M,\,n\leq N;\;m^u \equiv an^v\pmod q\},
$$
where, if $v$ is negative, then $n^{v}$ stands for ${\bar n}^{|v|}$. Moreover, $\bar n$ denotes the multiplicative inverse of $n$ modulo $q$.

It is not hard to see that one has the upper bound
\begin{equation}\label{nothard}
S_{u,v}(M,N,q,a)\ll \frac{MN}{q} +\min(M,N).
\end{equation}
In certain cases, this can even be improved by making use of the Weil bound for exponential sums over curves. For example, suppose $q$ is squarefree and $(u,v)=1$, $u\neq v$. Then we have the inequality
\begin{equation}\label{Weil}
S_{u,v}(M,N,q,a)\ll q^{\epsilon}\left( MNq^{-1} + Mq^{-\frac12} + Nq^{-\frac12} + q^{\frac12}\right).
\end{equation}
Unfortunately, when $M \asymp N \asymp q^{1/2}$, both \eqref{nothard} an \eqref{Weil} give the same bound
\begin{equation}\label{PV}
S_{u,v}(M,N,q,a)\ll q^{\frac12+\epsilon}.
\end{equation}
This is an important threshold when trying to improve \eqref{Hoo-sqf} or \eqref{3/2}. Indeed, one of the main achievements in \citep{heath1982least} is giving an upper bound for $S_{1,-2}(M,N,q,a)$ that improves on \eqref{PV} in the range where $M$ and $N$ are close to $q^{\frac12}$ in logarithmic scale.
The following lemma is a particular case of a result by Pierce \citep{pierce2005part} generalizing the main argument in \citep{heath1982least}. Both of these results are inspired by work of Burgess \citep{burgess1962characters}.
\begin{lem}\label{Pierce} (see \citep[Theorem 4]{pierce2005part})
We have, uniformly for $a\in (\Z/q\Z)^{\ast}$, $1\leq M\leq q^{3/4}$ and $1\leq N < q/2$ the inequality
$$
S(M,N,q,a)\ll M^{\frac{2}{3}}N^{\frac{1}{4}}q^{\varepsilon}.
$$
\end{lem}

We will use Lemma \ref{Pierce} with $(u,v)=(1,-2)$ and $(u,v)=(2,-1)$. For the first of these pairs, the work of Heath-Brown suffices and if we only had Lemma \ref{Pierce} for this value of $(u,v)$, we could prove a version of Theorem \ref{25/36} with the exponent $25/36$ replaced by $9/13$. Hence Corollary \ref{least} would be just a particular case of \citep[Theorem 2]{heath1982least}. It is thanks to the more powerful result from \citep{pierce2005part} and the simple symmetry relation
\begin{equation}\label{symmetry}
S_{u,v}(M,N,q,a) = S_{-v,-u}(N,M,q,a),
\end{equation}
that we can obtain the improved exponent 25/36.

\section{Initial steps}

Let $q$ be a squarefree number, let $a$ be coprime with $q$ and $X\geq q$. We consider $E=E(X,q,a)$ given by

$$
E:=\sum_{\substack{n\leq X\\n\equiv a \rmod q}}\mu^2(n) - \frac{1}{\varphi(q)}\sum_{\substack{n\leq X\\(n,q)=1}}\mu^2(n).
$$
Our goal is to prove that we have the inequality $E\ll X^{1-\delta}/q$ uniformly for $q\leq X^{\frac{25}{36}-\epsilon}$.

If $q\leq X^{1/2}$, then this already follows from \eqref{Hoo-sqf}. Therefore we may suppose
$$
q\geq X^{1/2}.
$$
We recall the classical identity
\begin{equation}\label{mu-decomp}
\mu^2(n)=\underset{\substack{n_1,n_2\geq 1\\n_1n_2^2=n}}{\sum\sum}\mu(n_2).
\end{equation}
This gives
$$
E=\sum_{\substack{n\leq X^{1/2}\\(n,q)=1}}\mu(n)\Delta(X/n^2,q,a{\bar n}^2),
$$
where for every $x\geq 1$, $q$ integer and $a\in \Z/q\Z$,
$$
\Delta(x,q,a):=\sum_{\substack{m\leq x\\m\equiv a \rmod q}}1 - \frac{1}{\varphi(q)}\sum_{\substack{m\leq x\\(m,q)=1}}1.
$$
It is clear that for any $x,q,a$, we have
$$
\Delta(x,q,a) \ll 1.
$$
Let $N_0$ be a parameter to be chosen optimally later and such that $1\leq N_0\leq X^{1/2}$. The previous inequality shows us that
\begin{equation}\label{N0-out}
E=\sum_{\substack{N_0<n\leq X^{1/2}\\(n,q)=1}}\mu(n)\Delta(X/n^2,q,a{\bar n}^2) + O(N_0q^{\epsilon}).
\end{equation}
Notice that
\begin{align*}
\frac{1}{\varphi(q)}\sum_{\substack{N_0<n\leq X^{1/2}\\(n,q)=1}}\mu(n)
\sum_{\substack{m\leq X/n^2\\(m,q)=1}}1\ll \frac{X^{1+\epsilon}}{N_0q},
\end{align*}
This and \eqref{N0-out} combined give
\begin{equation}\label{beforedyadic}
|E|\leq \sum_{\substack{N_0<n\leq X^{1/2}\\(n,q)=1}}\sum_{\substack{m\leq X/n^2\\m\equiv a{\bar n}^2\rmod q}}1 + O\left(X^{\epsilon}\left(N_0+\frac{X^{1+\epsilon}}{N_0q}\right)\right).
\end{equation}
\section{Division in dyadic boxes}
We now proceed by means of a dyadic decomposition. If we put
\begin{equation}\label{SV}
S(M,N,q,a)=\underset{\substack{m\sim M,\,n\sim N\\ mn^2\equiv a\rmod{q}}}{\sum\sum}1,
\end{equation}
we deduce from \eqref{beforedyadic} the upper bound
$$
E\ll (\log X)^2\cdot \sup_{M,N}S\left(M,N,q,a\right) + N_0 + \frac{X^{1+\epsilon}}{N_0 q},
$$
where the supremum is taken over all $M$ and $N$ such that
\begin{equation}\label{cond1}
M,N\geq 1,\, N_0\leq N\leq 2X^{1/2},\,MN^2\leq 8X.
\end{equation}
Let $M_0\geq 1$ be a parameter to be chosen optimally later. Suppose that $M\leq M_0$ and that $M,N$ satisfy the conditions \eqref{cond1}. Then, by the crude estimate 
$$
\sum_{\substack{n\sim N\\n \equiv \alpha\rmod{q}}}1\ll \frac{N}{q} + 1,
$$
we see that
\begin{align*}
S(M,N,q,a)\ll Mq^{\epsilon}\left(\frac{N}{q}+1\right)\\
\ll X^{\epsilon}\left(\frac{X}{N_0 q} + M_0\right).
\end{align*}
Thus
\begin{equation}\label{E-sup}
E\ll (\log X)^2\sup_{M,N}S\left(M,N,q,a\right) + X^{\epsilon}\left(M_0 + N_0 + \frac{X}{N_0 q}\right),
\end{equation}
where now the supremum is taken over all $M$ and $N$ satisfying
\begin{equation}\label{cond2}
M\geq M_0,\, N\geq N_0,\, MN^2\leq 8X.
\end{equation}
\section{Using Lemma \ref{Pierce}}
We notice that
\begin{equation}\label{bounded}
S(M,N,q,a)\leq S_{1,-2}(M,N,q,a).
\end{equation}
Suppose that $M_0$ and $N_0$ satisfy
\begin{equation}\label{toverify}
M_0> Xq^{-3/2}, N_0> X^{1/2}q^{-3/8}.    
\end{equation}
This readily implies that every $M,N$ satisfying \eqref{cond2} we have $1\leq M,N\leq q^{3/4}$. Lemma~\ref{Pierce}, \eqref{bounded} and \eqref{symmetry} now give the upper bound
\begin{equation}\label{min}
S(M,N,q,a)\ll q^{\epsilon}\min\left(M^{\frac{2}{3}}N^{\frac{1}{4}},M^{\frac{1}{4}}N^{\frac{2}{3}}\right)
\end{equation}
for every $M,N$ satisfying \eqref{cond2}.

It is not hard to see that for every $0<\alpha<1$, It follows from \eqref{min} that we have the inequality
\begin{equation}\label{alfalfa}
S(M,N,q,a)\ll q^{\epsilon}\min\left(M^{\frac{2}{3}}N^{\frac{1}{4}}\right)^{\alpha}\left(M^{\frac{1}{4}}N^{\frac{2}{3}}\right)^{1-\alpha}.
\end{equation}
Taking $\alpha=2/15$, we get
$$
S(M,N,q,a)\ll q^{\epsilon} \left(MN^2\right)^{\frac{11}{36}} \leq X^{\frac{11}{36}+\epsilon}.
$$
Now by \eqref{E-sup} we see that
\begin{equation}\label{almostthere}
E\ll  X^{\frac{11}{36}+\epsilon}+ M_0 + N_0 + N^{-1}_0Xq^{-1}.
\end{equation}
\section{Conclusion}
We make the choices
\begin{equation}\label{choices}
M_0=2\max(Xq^{-\frac{3}{2}},1),\,N_0=2 X^{\frac{1}{2}}q^{-\frac{3}{8}}.
\end{equation}
Note that these choices clearly satisfy \eqref{toverify}. We also notice that we have $1\leq M_0\leq X$ and $1\leq N_0\leq X^{1/2}$.
With the choices \eqref{choices}, the upper bound \eqref{almostthere} becomes.
$$
E\ll X^{\epsilon}\left( X^{\frac{11}{36}} + Xq^{-\frac32} +X^{\frac12}q^{-\frac38}\right).
$$
It is now straightforward to verify that for every $\epsilon>0$, there exists $\delta=\delta(\epsilon)>0$ such that whenever $q\leq X^{\frac{25}{36}-\epsilon}$, we have the inequality

$$
E\ll \frac{X^{1-\delta}}q.
$$
This concludes the proof of Theorem \ref{25/36}.
\bibliographystyle{amsplain}
\bibliography{references}
\end{document}